\documentclass[12 pt]{article}   

 \usepackage{amsmath}
 \usepackage{amssymb}
 \usepackage{amscd}
 \usepackage[mathscr]{eucal}
 \newtheorem{thm}{Theorem}

 \newtheorem{prop}[thm]{Proposition}

\newcommand{\R}{\mathbb{R}} 
 \newcommand{\Z}{\mathbb{Z}}
  
 \newcommand{\C}{\mathbb{C}}
\newcommand{\HH}{\mathbb{H}} 
\newcommand{\N}{\mathbb{N}}

 \renewcommand{\H}{\mathscr{H}}
 \renewcommand{\N}{\mathscr{N}}
 \newcommand{\U}{\mathscr{U}}  
 \newcommand{\A}{\mathscr{A}}
\newcommand{\B}{\mathscr{B}}

\begin{document}

\thispagestyle{empty}

\title{Intersection forms of  toric hyperk\"ahler  varieties}

\author{
 Tam\'as Hausel
\\ {\it  Department of Mathematics}
\\ {\it University of Texas at Austin}
\\ {\it Austin, TX 78712-1082, USA}
\\{\tt hausel@math.utexas.edu}
\\
\\Edward Swartz
\\ {\it Department of Mathematics}
\\ {\it Cornell University}
\\ {\it Ithaca NY 14853-4201, USA}\\
{\tt ebs@math.cornell.edu } }

\maketitle
\begin{abstract}
This note proves combinatorially that the intersection pairing on the middle dimensional compactly supported cohomology of a toric hyperk\"ahler variety is always definite, providing a large number of non-trivial $L^2$ harmonic forms for toric hyperk\"ahler metrics on these varieties. 
This is motivated by a result of Hitchin
about the definiteness of the pairing of $L^2$ harmonic forms on complete 
hyperk\"ahler manifolds of linear growth.
\end{abstract}

\section{Introduction}

A complete 
hyperk\"ahler manifold $M^{4n}$ has linear growth if for 
one of the K\"ahler forms 
$\omega=d\beta$ with $\beta$ of linear growth (see \cite{hitchin} for details).
In \cite{hitchin}[Theorem 4] Hitchin proved that on a complete 
hyperk\"ahler manifold $M^{4n}$ of linear growth non-trivial 
$L^2$ harmonic 
forms are of middle dimension and anti-self-dual (resp. self-dual) if $n$ is 
odd (resp. $n$ is even). It follows that the intersection form on the 
Hodge cohomology (the space of $L^2$ harmonic forms) on such manifolds 
is always definite (negative definite if $n$ is odd, positive definite if $n$ is even). 
Moreover, Hodge theory implies that  the image of middle dimensional 
compactly supported cohomology in ordinary cohomology
filters through  $L^2$ harmonic forms. 
We think of the natural intersection form on the image of compactly supported cohomology in 
ordinary cohomology as the intersection form, given by cup product and integration, 
on compactly supported cohomology modulo its 
null-space. 
We can thus conjecture that for a complete hyperk\"ahler manifold of linear growth and dimension $4n$ 
this intersection form is definite; in particular
this would imply that the signature of the manifold is non-positive if $n$ is odd and non-negative if $n$ is even. 
If true this  would be
a non-trivial topological obstruction for a manifold to carry a 
complete hyperk\"ahler metric of linear growth. 

This is known to 
hold in all the examples where the intersection form on middle dimensional 
compactly supported cohomology of a complete hyperk\"ahler manifold of linear 
growth has been calculated. These examples include Nakajima's deep result \cite{nakajima-quiver}[Corollary 11.2] that the 
intersection form on a quiver variety is always definite, which is deduced from
a representation theory result of Kac. In another example, Segal and Selby  \cite{segal-selby}
proved that the intersection form on the moduli space of $SU(2)$ magnetic monopoles on $\R^3$
is definite (fitting very nicely with string theory 
conjectures of Sen \cite{sen}
on the Hodge cohomology of such magnetic monopole moduli spaces). The final 
example is \cite{hausel}, where 
it was shown that the moduli space of rank $2$ Higgs bundles of fixed 
determinant of odd degree, another gauge theoretic example of a complete
hyperk\"ahler manifold of linear growth, has a trivial and so definite intersection form. 

A general theorem, \cite{hausel-etal}[Corollary 7],  
implies that when the metric of the complete hyperk\"ahler manifold of linear growth $M^{4n}$  is also of fibered boundary (or fibered cusp) type, then the 
intersection form  on middle dimensional compactly supported cohomology is 
always 
semi-definite; and consequently that
the intersection form on the image of compactly supported
cohomology in ordinary cohomology is definite. 
The examples of such manifolds include all
known gravitational instantons of finite topological type as well as smooth 
generic 
ALE toric hyperk\"ahler manifolds (cf. \cite{hausel-etal}[Section 7.2]). 
In fact, using our main Theorem~\ref{main} below, 
\cite{hausel-etal}[Corollary 11] 
calculates the full Hodge cohomology of such ALE toric hyperk\"ahler manifolds.

In this paper we study the intersection form of another family of 
complete hyperk\"ahler manifolds of linear growth. Bielawski and Dancer in \cite{bielawski-dancer} 
construct  toric hyperk\"ahler manifolds as hyperk\"ahler quotients of flat 
quaternionic 
space $\HH^n$ by a 
hyperk\"ahler action of a torus $T^d\subset T^n$. An account of the algebraic geometry 
of the underlying varieties from a combinatorial perspective 
were given in \cite{hausel-sturmfels}. We will
follow the notations and terminology of \cite{hausel-sturmfels}. 

The following is the main result of this note.

\begin{thm}
\label{main} Let $\theta\in {\mathbb N}\A$ be a smooth degree, so that the toric hyperk\"ahler variety $Y(A,\theta)$ of 
real dimension $4n-4d$ is smooth. 
Then the intersection form on $H^{2n-2d}_{cpt}(Y(A,\theta),\R)$ is always a 
definite form: it is positive definite if $n-d$ is
even, negative definite if $n-d$ is odd. 
\end{thm}

When the underlying hyperplane arrangement is co-graphic
the toric hyperk\"ahler variety is also a toric quiver variety (see 
\cite[Section 7]{hausel-sturmfels}). In these circumstances Nakajima's above mentioned result already proves the theorem. 

In the general case we proceed by considering the bounded complex of the 
affine hyperplane arrangement $\H^{bd}(\B,\psi)$ in $\R^{n-d}$ 
defined from the data $(A,\theta)$ (for details see \cite{hausel-sturmfels}). 
Then a basis for $H^{2n-2d}_{cpt}(Y(A,\theta),\R)$ is given by the compactly supported cohomology classes $\eta_{X_F}$
of middle dimensional projective subvarieties $X_F$ of $Y(A,\theta)$. 
Each $X_F$ is a toric variety associated to a top dimensional
bounded region $F$ in $\H^{bd}(\B,\psi)$.  We will show in the next section
that in this basis the intersection form is combinatorially given by:

\begin{thm} 
\label{intersection} \begin{equation}\int_{Y(A,\theta)} \eta_{X_1}\wedge \eta_{X_2} = (-1)^{\dim \overline{F_1}\cap \overline{F_2} }\left(\mbox{number of vertices of } \overline{F_1}\cap \overline{F_2}\right),\label{inter}\end{equation}
where $F_1$ and $F_2$ are two top dimensional bounded regions in $\H^{bd}(\B,\psi)$ and $X_1$ and $X_2$ are the corresponding projective toric varieties in $Y(A,\theta)$. Moreover the classes $\eta_{X_F}$, where $F$ runs through the top (i.e. $n-d$) 
dimensional
bounded regions in $\H^{bd}(\B,\psi)$ form a basis for the vector space 
$H^{2n-2d}_{cpt}(Y(A,\theta),\R)$.
\end{thm}

In the last section we then prove that this combinatorial 
intersection pairing, given purely in terms of the affine hyperplane arrangement $\H(\B,\psi),$ is indeed definite. We will
in fact construct a natural isomorphism of this pairing with the natural pairing on the $(n-d-1)$ 
cohomology of the independence complex $\N$ of the matroid $M(\B)$. 

\paragraph{\bf Acknowledgement.} The first author was 
partly supported by NSF grant DMS-0072675.  The second author was partly supported by  a VIGRE postdoc under NSF grant number 9983660
to Cornell University.

\section{Determining the intersection form} 

In this section we prove Theorem~\ref{intersection}. We will use the terminology and notation of \cite{hausel-sturmfels}. 
Let us fix $A$ and a smooth degree $\theta$, so that $\B$ is coloop-free. Thus  we are in the situation of \cite{hausel-sturmfels}[Proposition 6.7]. We can
assume this,  otherwise $H^{2n-2d}(Y,\R)=0$ and 
Theorem~\ref{intersection} holds automatically. 
For convenience we will write $Y$ for $Y(A,\theta)$ and $\H^{bd}$ for $\H^{bd}(\B,\psi)$.

First we note that the subvariety $X_F$ of $Y$ corresponding to a top dimensional bounded region $F$ in $ \H^{bd}$ is a Lagrangian subvariety with respect to the natural 
holomorphic symplectic structure $\omega_\C$ on
$Y$. One way to see this is to
use a circle action on $Y$ corresponding to the region $F$ as explained in \cite{harada-proudfoot}. By construction the holomorphic symplectic form
on $Y(A,\theta)$ is of homogeneity $1$ with respect to this circle action, meaning
that $\lambda^{*}(\omega_\C)=\lambda \omega_\C$, where $\lambda\in U(1)$. 
Moreover,  $X_F$ is the minimum of the associated moment map. Now it is clear from \cite{nakajima-book}[Proposition 7.1] that $X_F$ is a Lagrangian subvariety.

Now consider $F_1$ and $F_2$ two top dimensional bounded regions in $\H^{bd}$ and let $X_1$ and $X_2$ denote the corresponding projective
toric varieties in $Y$. Then the third equation in \cite{fulton-intersection}[Proposition 9.1.1] implies
that $$\int_Y \eta_{X_1}\wedge \eta_{X_2} = \int_{X_{12}} c(N_1)
\wedge c(T_{X_2})^{-1}\wedge c(T_{X_{12}}),$$ where $X_{12}$ is 
the projective toric variety 
in $Y$ corresponding to the region $\overline{F_1}\cap \overline{F_2}$, and $N_1$ denotes the normal
bundle of $X_1$ in $Y$. 

Since $X_1$ and $X_2$ are Lagrangian subvarieties, we see that $c(N_1)
\wedge c(T_{X_2})^{-1}\wedge c(T_{X_{12}}) = c(T^*(X_{12}))$ on $X_{12}$. Therefore,
 \begin{eqnarray*} \int_Y \eta_{X_1}\wedge \eta_{X_2}=\int_{X_{12}} c(T^*X_{12})&=&(-1)^{\dim(X_{12})} \chi(X_{12}) \\ &=& (-1)^{\dim(\overline{F_1}\cap \overline{F_2})}\left(\mbox{number of vertices of } \overline{F_1}\cap \overline{F_2})\right), \end{eqnarray*} where we used that for the Euler characteristic of toric variety we have 
$$\chi(X_{12})=(\mbox{number of vertices of } \overline{F_1}\cap\overline{F_2}).$$ 

To prove the last statement of our theorem list 
$F_1,F_2,\dots, F_r$ the top dimensional bounded regions of 
$\H^{bd}$ and the corresponding $X_1,\dots,X_r$ middle dimensional smooth projective subvarieties of $Y$. 
Then as in (35) of \cite{hausel-sturmfels}, we can find a basis 
$\alpha_1,\dots,\alpha_r$ for $H^{2n-2d}(Y,\R)$ which has the property that 
$\alpha_i\mid_{X_j}\neq 0$ if and only if $i=j$. The Poincar\'e dual basis for
$H_{cpt}^{2n-2d}(Y,\R)$ is then clearly $\eta_{X_1},\dots,\eta_{X_r}$. This 
completes the proof of Theorem~\ref{intersection}.

\section{Combinatorial intersection pairing}

In this section we prove that the combinatorial intersection pairing given by (\ref{inter}) is $(-1)^{n-d}$ times a positive definite pairing.  This will be done by using the classical nerve construction which we now review.

Let $\Gamma$ be a finite CW-complex and let $\Gamma_1,\dots,\Gamma_r$ be  subcomplexes of $\Gamma$ such that all  non-void intersections are contractible and the union of the $\Gamma_i$ cover $\Gamma.$  Let $\N$ be the {\it nerve} of this cover.  That is, $\N$ is the abstract simplicial complex with  vertices $v_1,\dots,v_r$ and whose simplices consist of all $[v_{i_1},\dots, v_{i_k}]$ such that $\Gamma_{i_1} \cap \dots \cap \Gamma_{i_k} \neq \emptyset.$  Let $Z$ be the subset of $\Gamma \times \N$ defined by

$$\{(x,z) : x \in \Gamma_{i_1} \cap \dots \cap \Gamma_{i_k} \text{ and } z \in [v_{i_1},\dots,v_{i_k}]\}.$$

\begin{thm} \cite{BS} \label{nerve} 
Let $\Gamma,\N$ and $Z$ be as above.  The projections $\pi_\Gamma$ and $\pi_\N$ of $Z$ onto $\Gamma$ and $\N$ 
respectively are homotopy equivalences.
\end{thm}

The bounded complex of any affine hyperplane arrangement is contractible (for a proof see Theorem 3.3 and Theorem 4.7 of \cite{hausel-sturmfels}, it also follows from \cite[Exercise 4.27 (a)]{BLSWZ}), hence $\H^{bd}(B,\psi)$ is contractible.  Let $\Gamma = \H^{bd}(B,\psi) - \bigcup F_j.$  Since the $F_j$ are disjoint, $\Gamma$ is homotopy equivalent to a wedge of $r$ spheres of dimension $n-d-1.$ 

Let $\{H_1, \dots, H_s \}$ be the affine hyperplanes in $\H.$ Define $\Gamma_i = \Gamma \cap H_i$ and $\U = \{\Gamma_1,\dots,\Gamma_s\}.$ 
Now, $\U$ covers $\Gamma$ and any non-void intersection of members of $\U$  is the bounded complex of an affine hyperplane arrangement, and hence is contractible.  Let $\N$ be the nerve of this cover.  We denote the vertices of $\N$ by $v_1, \dots, v_s.$  A subset of vertices of $\N$ is a simplex if and only if the corresponding hyperplanes have nonempty intersection and this holds if and only if the corresponding columns of $B$ are independent.  Since this only depends on the matroid $M(B), \ \N$ is known as the independence (or matroid) complex of the matroid $M(B).$  

Let $\sigma$ be a $k$-dimensional cell of $\Gamma$ and let $(H_{\sigma_1}, \dots, H_{\sigma_{n-d-k}})$ be the ordered set of hyperplanes containing $\sigma.$  For each hyperplane $H_{\sigma_i}$ choose a normal $\eta_{\sigma_i}.$ The ordered set of normals $(\eta_{\sigma_1},\dots,\eta_{\sigma_n-d-k})$ define an orientation for $\sigma$ as follows.  If $x \in \sigma$ and $W=(w_1,\dots,w_k)$ is an ordered basis for the tangent space of $\sigma$ at $x,$ then $W$ is positively oriented if and only if $(\eta_{\sigma_1},\dots,\eta_{\sigma_n-d-k},w_1,\dots,w_k)$ is a positively oriented basis for the tangent space of $\R^{n-d}$ at $x.$

Let $F$ be an $(n-d)$-dimensional cell of $\H^{bd}(B,\psi).$  Choose an inward pointing normal for each hyperplane incident to $F.$ Define $[F]$ to be the cycle in $H_\star(\Gamma,\Z)$ (CW-homology) given by $\sum [\sigma],$ where the sum is taken over all cells $\sigma$ of dimension $n-d-1$ on the boundary of  $F.$  For each vertex $v$ on the boundary of $F$ let  $\Psi(v)$ to be the oriented simplex $[v_{i_1},\dots,v_{i_{n-d}}]$ in $\N,$ where $H_{i_1},\dots,H_{i_{n-d}}$ are the hyperplanes which contain $v,$ and the orientation of the simplex is determined by the orientation of the vertex as a zero-cell of $\Gamma.$  Now define

$$\Psi[F] = \sum_{v \in \partial F} \Psi(v).$$

\noindent In order to see that $\Psi[F]$ is a cycle,  note that the facets  in $\partial \Psi[F]$ correspond to the one-cells in $\partial F$ and these occur in oppositely oriented pairs, one for each end point of the one-cell.  Since $\{[F_1],\dots,[F_r]\}$ is a basis for $H_\star(\Gamma;\Z),$ we can extend $\Psi$ linearly to a map $\Psi_\star:H_\star(\Gamma;\Z) \to H_\star(\N;\Z).$

\begin{prop} \label{combinatorial definiteness}
  $\Psi_\star$ is an isomorphism.  
\end{prop}  

{\em Proof:} Let $Z = \{(x,z) \in \Gamma \times \N : x \in \Gamma_{i_1} \cap \dots \cap \Gamma_{i_k} \text{ and } z \in [v_{i_1},\dots,v_{i_k}]\}.$ Fix $F.$ By Theorem \ref{nerve} it is sufficient to find a cycle $[\zeta] \in H_\star(Z;\Z)$ such that $(\pi_\Gamma)_\star[\zeta] = [F]$ and $(\pi_{\N})_\star[\zeta] = \Psi[F].$  

In order to define $\zeta$ we introduce the following notation.  Let $\Delta = [v_{i_1},\dots,v_{i_k}]$ be a (non-empty) simplex of $\N.$  Denote by $H_\Delta$ the corresponding  $(n-d-k)$-cell in $\Gamma$  with orientation given by $(\eta_{i_1},\dots,\eta_{i_k}).$ Furthermore,  let $(H_\Delta,v_j)$ be the cell $H_\Delta \cap H_j$ with orientation given by $(\eta_{i_1},\dots,\eta_{i_k},\eta_j).$ If $v_j \in \Delta$ or  $\Delta \cup v_j$ is not a simplex of $\N,$ then $(H_\Delta, v_j)$ is empty and the chain $[(H_\Delta,v_j)]$ equals zero. Define a cellular chain $\zeta$ in $C_\star(Z)$ by

$$\zeta = \sum_{H_\Delta \subseteq \partial c} (-1)^{|\Delta|} [H_\Delta \times \Delta].$$

\noindent Recall that for products of CW-complexes the boundary map is given by $\partial(\Omega \times \Psi) = \partial \Omega \times \Psi + (-1)^p \Omega \times \partial \Psi,$ where $\Omega$ is a $p$-cell and $\Psi$ is a $q$-cell.  Therefore,

$$\begin{array}{lcl}
 \vspace{2 mm}\partial^2(\zeta) & = &  \displaystyle\sum_j \partial([(H_\Delta,v_j) \times \Delta] + (-1)^{|\Delta|}[H_\Delta \times \partial \Delta]) \\
\vspace{2 mm}\ &= & \displaystyle\sum_{j,l} ([(H_\Delta,v_j,v_l) \times \Delta] + (-1)^{|\Delta + 1|}[(H_\Delta,v_j) \times \partial \Delta] \\ 
\ & \  & \ \ \ \ \ + \ (-1)^{|\Delta|} [(H_\Delta,v_j) \times \partial \Delta] + [H_\Delta \times \partial^2 \Delta]) \\
\ & = & 0.
\end{array} $$

\noindent Thus $\zeta$ is a cycle.  It is easy to see that $(\pi_\Gamma)_\star[\zeta] = [F]$ and $(\pi_{\N})_\star[\zeta] = \Psi[F].  \hfill \square$

The chain complex $C_\star(\N)$ has an inner product structure given by declaring that the set of chains $\{[\Delta]: \Delta \mbox{ a simplex of } \N\}$ is an orthonormal basis. Since $\N$ is a simplicial complex of dimension $n-d-1$ and is homotopy equivalent to  a wedge of $(n-d-1)$-dimensional spheres, $H_{n-d-1}(\N)$ is a subspace of $C_{n-d-1}(\N)$ and inherits a positive definite inner product.  

For convenience we recall the combinatorial pairing introduced in (\ref{inter}) in the present notation. Let $V(B)$ be the vector space whose basis is $F_1, \dots, F_r.$ Set $\sigma_{ij} = \overline{F_i} \cap \overline{F_j}.$  Define a pairing $\Phi(F_i,F_j) = (-1)^{\dim \sigma_{ij}} | \{ \mbox{vertices } v: v \in \sigma_{ij} \}|.$

\begin{prop}
\label{identification}
  $\Phi(F_i,F_j) = (-1)^{n-d} <\Psi(F_i),\Psi(F_j)>.$
\end{prop}

{\em Proof:}
  It is evident that $<\Psi(F_i),\Psi(F_j)>$ counts the number of vertices in $\sigma_{ij}$ weighted with $+1$ if the orientations induced on the vertex by the inward pointing normals of $F_i$ and $F_j$ are the same, $-1$ if they are different.  The orientation of a vertex with respect to the inward pointing normals for $F_j$ is obtained from the orientation of the vertex with respect to $F_i$ by reversing the direction of $n-d-\dim \sigma_{ij}$ of the normals.  Hence $<\Psi(F_i),\Psi(F_j)> = (-1)^{n-d-\dim \sigma_{ij}} | \{ \mbox{vertices } v: v \in \sigma_{ij} \}|. 
 \hfill \square$

\paragraph{\bf Remark} 1. The map $\Psi_\star$
identifies $H_*(\N,Z)$ with $H_*(\Gamma;Z)$, which in turn could naturally be identified with the middle dimensional compactly supported cohomology $H^{2n-2d}_{cpt}(Y,\Z)$. Moreover by Theorem~\ref{intersection} and Proposition~\ref{identification} these identifications also preserve the appropriate inner products
on these spaces, so as a byproduct we get Theorem~\ref{main}. 

  2. Through the above identifications 
the map $\Psi_\star$ is defining a flat connection on the
middle dimensional 
compactly supported cohomology of toric hyperk\"ahler varieties $Y(A,\theta)$ 
as  $A$ is fixed and $\theta$ varies. 
We can think of this as a combinatorial version of 
the Gauss-Manin connection obtained from the hyperk\"ahler quotient 
construction. 

3. While $\Psi_\star$ makes sense for arbitrary hyperplane arrangements, Proposition~\ref{combinatorial definiteness} may not hold if the arrangement is not generic.  In the arrangement pictured in figure 1 
 $$\Phi(F_1+F_2-F_3-F_4, F_1+F_2-F_3-F_4) <0.$$

\vspace{1in}

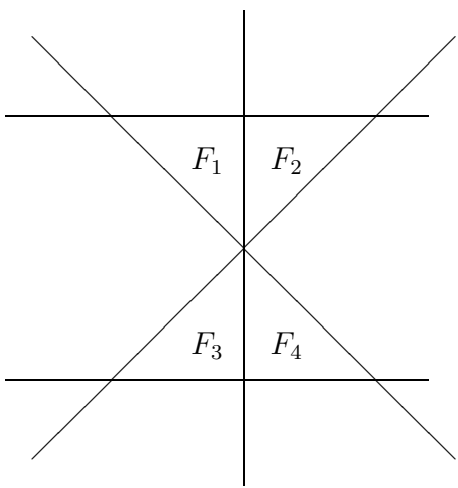
\begin{figure}[here] 

\begin{picture}(2500,100)(-250,-100)  
 \put(-90,50){\line(1,0){160}}
  \put(-90,-50){\line(1,0){160}}
  \put(-80,80){\line(1,-1){160}}
  \put(0,90){\line(0,-1){180}}
  \put(-80,-80){\line(1,1){160}}
  \put(-20,30){$F_1$}
   \put (10,30){$F_2$}
  \put(-20,-40){$F_3$}
  \put(10,-40){$F_4$}
 \end{picture} 

\caption{A nongeneric arrangement}
\end{figure}

\end{document}